\documentclass[10pt]{article}
\usepackage{amsfonts}
\usepackage{amsmath}
\usepackage{amscd,amssymb,amsthm}
\usepackage{graphics}

\newtheorem{theorem}{\sc Theorem}[section]

\newtheorem{prop}{\sc Proposition}[section]
\newtheorem{remark}{\sc Remark  }[section]

\newcommand{\bce}{\begin{center}}
\newcommand{\ece}{\end{center}}

\begin{document}

\title{Analysis of a $3D$ chaotic system  }

\date{}
\author{Gheorghe Tigan \thanks{Department of Mathematics, "Politehnica" University of
Timisoara, P-ta Victoriei, Nr.2, 300006, Timisoara, Timis,
Romania, email: gheorghe.tigan@mat.upt.ro }\ ,
Dumitru Opri\c{s} \thanks{%
 Department of Mathematics, West University of Timisoara, B-dul V. Parvan, Nr. 4, 300223, Timisoara, Timis,
 Romania, email: opris@math.uvt.ro } } \maketitle

\begin{abstract}
A $3D$ nonlinear chaotic system, called the $T$ system, is
analyzed in this paper. Horseshoe chaos is investigated via the
heteroclinic Shilnikov method constructing a heteroclinic
connections between the saddle equilibrium points of the system.
Partially numerical computations are carried out to support the
analytical results.\end{abstract}

 \section{Introduction}

 One of the first examples of continuous dynamical system of
 dimension three, whose numerical simulations display the
 property of sensitivity to initial conditions, is the Lorenz
 system \cite{lore1}:

\begin{equation}\label{1} \dot{x}=a(y-x), \ \dot{y}=cx-y-xz,
\ \dot{z}=-bz+xy,
\end{equation}
\noindent where $a,b,c$ are positive real parameters.  Physical
realization of the Lorenz system is the Rayleigh-Benard
experiment. The system was derived from the hydrodynamical
Navier-Stokes equations to create a dynamical model for
meteorology.
\par Trying to transform  the Lorenz system from a stable state to
a chaotic one (concept known as anticontrol of chaos or
chaotification, \cite{keng} ), Chen in \cite{chenueta1} introduces
a 3D polynomial system (of type Lorenz) known as the Chen system.
Vanecek and Celikovsky  in \cite{Vanecek}, characterize a
generalized Lorenz system by a condition on its linear part matrix
$J=(a_{ij})$: $\,\,a_{12}a_{21}>0$. L\"{u} in \cite{JinChe1}
observes that while the classical Lorenz system belongs to this
class, the Chen system satisfies a \emph{dual} condition
$a_{12}a_{21}<0$, (systems that satisfy such a condition are
called dual systems to the Lorenz system) and introduces a new
system that bridges the gap between the Chen system and Lorenz
system, satisfying $a_{12}a_{21}=0$. A nonlinear system arising
from a nuclear spin generator and compared with the Lorenz system
is studied in \cite{Sach}. Nonlinear dynamics is met in many areas
from Ecology \cite{ecology1} to Physics \cite{nasc1}.
\par From the point of view of the potential applications,
systems with sensitivity to the initial conditions can be used in
secure communications \cite{alvarez}, \cite{sprot}. Of the
pioneering papers which proposed to use the chaotic systems in
communications are the papers of Pecora and Carrol \cite{pec1},
\cite{pec2}. Consequently, an appropriate chaotic system can be
chosen from a catalogue of chaotic systems to optimize some
desirable factors, idea suggested in \cite{sprot}.

\par These some ideas led us to study a new 3D polynomial
differential system given by:
\begin{equation}\label{2}
\dot{x}=a(y-x), \ \dot{y}=(c-a)x-axz, \ \dot{z}=-bz+xy,
\end{equation}
\noindent with $a,b,c$ real parameters and $a\neq 0$. Call it the
$T$ system. Some results regarding the $T$ system are already
presented in \cite{tig1} and \cite{tig4}. Compared with the L\"{u}
system introduced in \cite{JinChe1},  the system $T$ allows a
larger possibility in choosing the parameters of the system and,
therefore, it displays a more complex dynamics.
\par The paper is organized as follows. In the second Section we
give some details of the equilibria. In Section 3, we pay
attention to the chaos in the system $T$ and, using a method of
type Shilnikov, we prove that the system displays "horseshoe"
chaos. This implies that the system possesses a strange chaotic
attractor. Not any strange dynamics is chaotical as it is pointed
out, for example, in \cite{nasc2}, where is presented a scenario
for a 2-dimensional dynamic system to possesses strange nonchaotic
behavior.

\section{The equilibrium points of the system T}

Straightforward computations lead us to the following result:

\begin{prop} If $\frac{b}{a}(c-a)>0$, the system $T$
possesses three equilibrium isolated points:\\
$O(0,0,0),E_{1}(\sqrt{\frac{b}{a}(c-a)},\sqrt{\frac{b}{a}(c-a)},\frac{c-a}{a}%
),E_{2}(-\sqrt{\frac{b}{a}(c-a)},-\sqrt{\frac{b}{a}(c-a)},\frac{c-a}{a})$, \\
\\
and for $b\neq 0,\frac{b}{a}(c-a)\leq 0 $ it has only one isolated
equilibrium point, $O(0,0,0).$
\end{prop}
As reported in \cite{tig1}, we have:
\begin{theorem}\label{t1} For $b\neq 0$ the following statements are true:\\
\\
\noindent a) If  $\left( a>0, b>0, c\leq a\right)$, then
$O(0,0,0)$ is asymptotically stable ,\\
\noindent b) If $\left( b<0\right) $ \ or \ $\left( a<0\right) $ \
or \ $\left( a>0 , c>a\right) $, then $O\left( 0,0,0\right) $ \ is
unstable.\end{theorem}

\begin{theorem}\label{t2} If $(a+b>0, ab(c-a)>0, b(2a^{2}+bc-ac)>0),$ the equilibrium points\\
 $E_{1,2}(\pm\sqrt{\frac{b}{a}(c-a)},\pm\sqrt{\frac{b}{a}(c-a)},\frac{c-a}{a})$ are asymptotically
stable.
\end{theorem}

\section{Chaos in the $T$ system  \emph{via} the heteroclinic Shilnikov method}

\par Computing the Lyapunov exponents with the software
\emph{Dynamics} \cite{nusse} for the parameter vector
$(a,b,c)=(2.1, 0.6, 30)$ and the initial conditions
$(0.1,-0.3,0.2)$, we get that $\lambda_1=0.37>0$, $
\lambda_2=0.00$ and $\lambda_3=-3.07$. So the system $T$ displays
chaotic characteristics.

\par In the following, using the Shilnikov heteroclinic method,
we show that the system $T$ presents chaos of \emph{horseshoe}
type. We will employ the following result:

\begin{theorem}(\cite{Silva, chen07}) If a $3D$ given system $\dot{x}=F(x)$ has two equilibrium points
$E_{1},$\ $E_{2}$, of type saddle-focus, i.e. the eigenvalues of
the Jacobi matrix associated to the system in these points are
$\gamma_{k}\in\mathbb{R}$ and $ \alpha_{k}\pm
i\beta_{k}\in\mathbb{C}$, $k=1,2,$ such that

\begin{equation}\label{horse2}
\alpha_{1}\alpha_{2}>0 \ \ or \ \ \gamma_{1}\gamma_{2}>0
\end{equation}
and (the Shilnikov inequality)
\begin{equation}\label{horse1}
|\gamma_{k}|>|\alpha_{k}|, k=1,2, \end{equation}

\noindent and if the system has a heteroclinic orbit connecting
the equilibrium points $E_{1}$ and $E_{2}$, then the Poincar\'{e}
map defined on a transversal section of the flow in a neighborhood
of the heteroclinic orbit presents chaos of horseshoe type.
\end{theorem}

We will show that our system $T$ given by
\begin{equation}\label{horse8}
\dot{x}=a(y-x), \ \dot{y}=(c-a)x-axz, \ \dot{z}=-bz+xy,
\end{equation}
fulfills the conditions from the above theorem.  The equilibrium
points of the system $T$ are:
$E_{1,2}(\pm\sqrt{\frac{b}{a}(c-a)},\pm\sqrt{\frac{b}{a}(c-a)},\frac{c-a}{a})$
for $\frac{b }{a}(c-a)>0$. \noindent The characteristic polynomial
associated to the Jacobi matrix of the system $T$ in these points
is:
\begin{equation}\label{ste1} f\left( \lambda \right)
=\lambda ^{3}+\lambda ^{2}(a+b)+bc\lambda +2ab(c-a)=0
\end{equation}
\noindent Denoting $\lambda=\mu-(a+b)/3$, (\ref{ste1}) leads to:

\begin{equation}\label{horse4}
\mu ^{3}+p\mu +q=0,
\end{equation}

\noindent where
\begin{equation}\label{horse5}
p=bc-\frac{1}{3}a^{2}-\frac{2}{3}ab-\frac{1}{3}b^{2}
\end{equation}
and
\begin{equation}\label{horse6}
q=\frac{2}{27}a^{3}-\frac{16}{9}a^{2}b+\frac{2}{9}\allowbreak ab^{2}+\frac{2%
}{27}b^{3}+\frac{5}{3}bca-\frac{1}{3}b^{2}c.
\end{equation}

Denote $\Delta =\left( \frac{q}{2}\right) ^{2}+\left(
\frac{p}{3}\right) ^{3}. $ Then, if $\Delta >0$, Eq.
(\ref{horse4}) has a negative solution, $\alpha_{1}$, together
with a pair of complex solutions,
$\alpha_{2}\pm i\alpha_{3}$, where \\

$$\alpha _{1}=\sqrt[3]{-\frac{q}{2}+\sqrt{\Delta }}+\sqrt[3]{-\frac{q}{2}-%
\sqrt{\Delta }},$$

$$\alpha_{2}=\frac{-1}{2}\left( \sqrt[3]{-\frac{q}{2}+\sqrt{\Delta }}+\sqrt[3]%
{-\frac{q}{2}-\sqrt{\Delta }}\right) ,$$

$$\alpha_{3}=\frac{\sqrt{3}}{2}\left( \sqrt[3]{-\frac{q}{2}+\sqrt{\Delta }}-%
\sqrt[3]{-\frac{q}{2}-\sqrt{\Delta }}\right) .$$ \\
So, if $\Delta >0$, the three roots of Eq. (\ref{ste1}) are:

\begin{equation}\label{horse7}
\lambda_{1}=-\frac{a+b}{3}+\alpha_{1}, \
\lambda_{2}=-\frac{a+b}{3}+\alpha_{2}+i\alpha_{3}, \
\lambda_{3}=-\frac{a+b}{3}+\alpha_{2}-i\alpha_{3}, \end{equation}
\noindent with $\lambda_{1}<0$.
 From the first two Eqs. of (\ref{horse8}) one gets

$y=x+\frac{\dot{x}}{a}$, \\

$z=-\frac{\dot{y}-\left( c-a\right) x}{ax}=-\frac{\ddot{x}+a\dot{x}}{a^{2}x}+%
\frac{c-a}{a}$, and using the last Eq. of (\ref{horse8}) we get \\

$\frac{d}{dt}\left( \frac{\ddot{x}+a\dot{x}}{x}\right) +b\frac{\ddot{x}+a%
\dot{x}}{x}+a^{2}x^{2}+ax\dot{x}-ab\left( c-a\right) =0,$
\\
or, equivalently
\\
\begin{equation}\label{horse9}
x\dddot{x}+\left( a+b\right) x\ddot{x}-\dot{x}\ddot{x}-a\dot{x}^{2}+abx%
\dot{x}+a^{2}x^{4}+ax^{3}\dot{x}-ab\left( c-a\right) x^{2}=0.
\end{equation}
We observe that if $x\left( t\right)$ is determined, then we can
find $y\left( t\right), z\left( t\right)$. Therefore, we have to
find a function $\phi \left( t\right)$ such that $x(t)=\phi(t)$ to
satisfy Eq. (\ref{horse9}) and $\phi(t)\rightarrow
-\sqrt{\frac{b}{a}(c-a)}$ for $t\rightarrow+\infty$,
$\phi(t)\rightarrow \sqrt{\frac{b}{a}(c-a)}$ for
$t\rightarrow-\infty$\ or,inversely $\phi(t)\rightarrow
\sqrt{\frac{b}{a}(c-a)}$ for $t\rightarrow+\infty$,
$\phi(t)\rightarrow
-\sqrt{\frac{b}{a}(c-a)}$ for $t\rightarrow-\infty$. \\
\noindent We can define, without loosing the generality, that the
direction from $E_{1}$ to $E_{2}$ corresponds to
$t\rightarrow\infty$, and from $E_{2}$ to $E_{1}$ corresponds to
$t\rightarrow-\infty$.

\par Consider  $\phi(t)$ of the following form $\phi \left( t\right)
=-x_{0}+\sum\limits_{n=1}^{\infty }a_{n}e^{n\alpha t},$ with
$\alpha $ real number and\\ $x_{0}=\sqrt{\frac{b}{a}\left(
c-a\right) }$.

\par Identifying the coefficients of $e^{n\alpha t}$ in
(\ref{horse9}) we get:
\begin{equation}\label{horse10}
x_{0}\left[ \alpha ^{3}+\alpha ^{2}(a+b)+bc\alpha +2ab(c-a)\right]
a_{1}=0, \end{equation}

\begin{equation}\label{horse11}
a_{2}=\frac{1}{f\left( 2\alpha \right) x_{0}}\left( \alpha
^{2}b+ab\alpha -ab(c-a)+6a^{2}x_{0}^{2}+3\alpha ax_{0}^{2}\right)
a_{1}^{2}, \end{equation}

\begin{equation}\label{horse12}
\begin{array}{cc}
  a_{3}=-\frac{1}{f\left( 3\alpha \right) x_{0}}\left(
4a^{2}x_{0}a_{1}^{3}+3a\alpha x_{0}a_{1}^{3}\right) + \\
\\
+ \frac{1}{f\left( 3\alpha \right) x_{0}}\left[ 3\alpha
^{3}+\alpha ^{2}a+5\alpha ^{2}b+3ab\alpha
-2abc+2a^{2}b+12a^{2}x_{0}^{2}+9\alpha ax_{0}^{2}\right]
a_{1}a_{2},
\end{array}%
\end{equation}

and for $n\geq 4$

\begin{equation}\label{horse13}
\begin{array}{cc}
  a_{n}=\frac{1}{f\left( n\alpha \right) x_{0}}\left(
a^{2}a_{ijpq}+ab_{ijpq}\right)+ \frac{1}{f\left( n\alpha \right)
x_{0}}\sum\limits_{i+j=n}\alpha ^{3}\left(
j^{3}-ij^{2}\right)a_{i}a_{j}+ \\
\\
+\frac{1}{f\left( n\alpha \right) x_{0}}
\sum\limits_{i+j=n}\left[\alpha ^{2}\left( \left( a+b\right)
j^{2}-ija\right) +ab\alpha j-ab(c-a)+6a^{2}x_{0}^{2}+3\alpha ajx_{0}^{2}%
\right] a_{i}a_{j},
\end{array} \end{equation}

\noindent where
\begin{equation}\label{horse14}
f\left( n\alpha \right) =\left( n\alpha \right) ^{3}+\left(
n\alpha \right) ^{2}(a+b)+bcn\alpha +2ab(c-a),\end{equation}

\begin{equation}\label{horse15}
a_{ijpq}=-4x_{0}\sum\limits_{i+j+p=n}a_{i}a_{j}a_{p}+\sum%
\limits_{i+j+p+q=n}a_{i}a_{j}a_{p}a_{q}, \end{equation}

\begin{equation}\label{horse16}
b_{ijpq}=-3\alpha
x_{0}\sum\limits_{i+j+p=n}pa_{i}a_{j}a_{p}+\alpha
\sum\limits_{i+j+p+q=n}pa_{i}a_{j}a_{p}a_{q},\end{equation} with
$i,j,p,q\geq 1.$

\bigskip

In order to can be implemented on the computer, we put the sums
(\ref{horse15}) and (\ref{horse16}) in the following forms:

\begin{equation}\label{horse17}
a_{ijpq}=-4x_{0}\sum\limits_{i=1}^{n-2}a_{i}\left(
\sum\limits_{j=1}^{n-i-1}a_{j}a_{n-i-j}\right)
+\sum\limits_{i=1}^{n-3}a_{i}\left(
\sum\limits_{j=1}^{n-i-1}a_{j}\left(
\sum\limits_{p=1}^{n-i-j-1}a_{p}a_{n-i-j-p}\right) \allowbreak
\allowbreak \right),\end{equation}

\begin{equation}\label{horse18}
\begin{array}{cc}
  b_{ijpq}=-3\alpha x_{0}\sum\limits_{i=1}^{n-2}a_{i}\left(
\sum\limits_{j=1}^{n-i-1}\left( n-i-j\right) a_{j}a_{n-i-j}\right)+ \\
\\
+ \alpha \sum\limits_{i=1}^{n-3}a_{i}\left(
\sum\limits_{j=1}^{n-i-1}a_{j}\left( \sum\limits_{p=1}^{n-i-j-1}p
a_{p}a_{n-i-j-p}\right) \right).
\end{array}
\end{equation}

Assume $a_{1}\neq 0$. If not, we observe inductively that all
coefficients $a_{n}, n\geq2$, are zero. Hence, from
(\ref{horse10}), we get
\begin{equation}\label{horse19}
\alpha^{3}+\alpha^{2}(a+b)+bc\alpha +2ab(c-a)=0
\end{equation} that is $\alpha$ is the negative root of
characteristic polynomial (\ref{ste1}). Because Eq. (\ref{ste1})
has a single negative solution for $\Delta>0$, we get that
\begin{equation}\label{horse20}
f\left( n\alpha \right) =\left( n\alpha \right) ^{3}+\left(
n\alpha \right) ^{2}(a+b)+bcn\alpha +2ab(c-a)\neq 0,
n>1.\end{equation} Consequently, the coefficients $a_{n}$ are
completely determined by $a,b,c,\alpha$ and $a_1$, and they are of
the following form:
\begin{equation}\label{horse21}
a_{n}=g(n)a_{1}^{n}, n>1
\end{equation} where the terms $g(n)$ are known functions. So the corresponding
branch of the heteroclinic orbit for $t>0$ is determined. In a
similar manner, consider the case $t<0$ and $\phi \left( t\right)
=x_{0}+\sum\limits_{n=1}^{\infty }a_{n}e^{-n\beta t},$ with
$\beta$ real, and find $\beta=\alpha$ and $b_{n}=-a_{n}, n>0$.
Consequently, the heteroclinic orbit $\phi(t)$ is given by:

\begin{equation}\label{horse22}
\phi \left( t\right) = \left\{
\begin{array}{ccc}
  -x_{0}+\sum\limits_{n=1}^{\infty }a_{n}e^{n\alpha t} & for &
  t>0 \\
  \\
  0 & for & t=0 \\
  \\
  x_{0}-\sum\limits_{n=1}^{\infty }a_{n}e^{-n\alpha t} & for & t<0 \\
\end{array} \right.
\end{equation}
Following a method presented in \cite{chen07} one can show that
$\phi(t)$ is uniformly convergent. So we have the following result
that characterizes the chaos in the system $T$:

\begin{theorem} If $\Delta>0$ and
$\alpha_{1}+\alpha_{2}<-\frac{2(a+b)}{3}$, then the system $T$ has
a heteroclinic orbit given by (\ref{horse22}), which connects the
equilibrium points $E_{1},E_{2}$, so the chaos is of horseshoe
type.
\end{theorem}

Let us partially illustrate the particular case $a=2.1, \ b=0.6, \
c=30.$ Then $\lambda_{1} =\alpha=-3.429$, $\lambda_{2}
=0.364-4.513i$, $\lambda_{3} =0.364+4.513i$, $\Delta = 911.69>0$,
$x_{0}= 2.8234$. Observe that the equilibrium points $E_{1},E_{2}$
are saddle-focus. Imposing that $\phi(t)$ to be at least continue,
from $\phi(0_{-})=\phi(0_{+})=0$, we find the equation
\begin{equation}\label{horse23}
-x_{0}+a_{1}+a_{2}+...+a_{n}=0, \ n>1.
\end{equation}
Solving this equation, we observe that the first coefficient is
$a_{1}=3.051$ for any $n>10$.

\begin{remark} Numerical series (\ref{horse22}) which describes
the heteroclinic orbit between $E_{1}$ and $E_{2}$, is rapidly
convergent. For example, considering the first ten terms of the
series, for $t=10$ we get
$\phi(t)=-x_{0}+\sum\limits_{n=1}^{10}a_{n}e^{n\alpha t}=-2.8234$,
for $t=-10$, $\phi(t)=x_{0}-\sum\limits_{n=1}^{10}a_{n}e^{-n\alpha
t}=2.8234$ and for $t=0$,
$\phi(t)=-x_{0}+\sum\limits_{n=1}^{10}a_{n}e^{n\alpha t}=0.0000$.
So the equilibrium points $E_{1}$, $E_{2}$ and $O$ are found with
an approximation of four exact decimals because
$x_{0}=\sqrt{\frac{b}{a}(c-a)}=2.8234$.\end{remark}

\section{Conclusions}
In this paper we investigated a new chaotic system. For some
parameter vectors, the system presents the property of dependence
to initial conditions which is a necessary condition for a system
to be chaotic. The system possesses three equilibrium points, the
origin $O(0,0,0)$ and another two points $E_{1,2}$. The chaos of
horseshoe type was analyzed using the heteroclinic Shilnikov
method, constructing a heteroclinic connection between the
saddle-focus equilibrium points $E_1$ and $E_2$.

\section{Acknowledgements}

The work of G.T. was supported through a research fellowship
offered by the Romanian Government by CNBSS and partially through
a European Community Marie Curie Fellowship, in the framework of
the CTS, contract number HPMT-CT-2001-00278.\\
The authors are grateful to Prof. M.S. El Naschie for his useful
suggestions.

\end{document}